\documentclass[preprint,number,sort&compress]{elsarticle}
\usepackage[utf8]{inputenc}
\usepackage{microtype} 

\usepackage{amsmath,amsfonts,amssymb,amsthm}
\usepackage{mathtools}
\usepackage{hyperref} 
\hypersetup{hidelinks}
\usepackage[nameinlink]{cleveref}
\usepackage{caption}
\usepackage{subcaption}

\hypersetup{pdfauthor=Jens Jäschke et al.}

\makeatletter
\DeclareFontFamily{OMX}{MnSymbolE}{}
\DeclareSymbolFont{MnLargeSymbols}{OMX}{MnSymbolE}{m}{n}
\SetSymbolFont{MnLargeSymbols}{bold}{OMX}{MnSymbolE}{b}{n}
\DeclareFontShape{OMX}{MnSymbolE}{m}{n}{
    <-6>  MnSymbolE5
   <6-7>  MnSymbolE6
   <7-8>  MnSymbolE7
   <8-9>  MnSymbolE8
   <9-10> MnSymbolE9
  <10-12> MnSymbolE10
  <12->   MnSymbolE12
}{}
\DeclareFontShape{OMX}{MnSymbolE}{b}{n}{
    <-6>  MnSymbolE-Bold5
   <6-7>  MnSymbolE-Bold6
   <7-8>  MnSymbolE-Bold7
   <8-9>  MnSymbolE-Bold8
   <9-10> MnSymbolE-Bold9
  <10-12> MnSymbolE-Bold10
  <12->   MnSymbolE-Bold12
}{}

\let\llangle\@undefined
\let\rrangle\@undefined
\DeclareMathDelimiter{\llangle}{\mathopen}%
                     {MnLargeSymbols}{'164}{MnLargeSymbols}{'164}
\DeclareMathDelimiter{\rrangle}{\mathclose}%
                     {MnLargeSymbols}{'171}{MnLargeSymbols}{'171}
\makeatother

\DeclarePairedDelimiter{\abs}{\vert}{\vert}
\DeclarePairedDelimiter{\norm}{\Vert}{\Vert}
\DeclarePairedDelimiterX{\scp}[2]{\langle}{\rangle}{#1,#2} 
\DeclarePairedDelimiterX{\sdp}[2]{\llangle}{\rrangle}{#1,#2} 

\newcommand{\pp}[2]{\frac{\partial #1}{\partial #2}} 
\newcommand{\dd}[2]{\frac{\dif{#1}}{\dif{#2}}} 
\newcommand{\Reals}{\mathbb{R}}
\newcommand*\dif{\mathop{}\!\mathrm{d}}

\DeclareMathOperator{\grad}{grad}
\DeclareMathOperator{\divergence}{div}

\theoremstyle{plain}
\newtheorem{theorem}{Theorem}[section]

\theoremstyle{definition}
\newtheorem{definition}[theorem]{Definition}

\theoremstyle{remark}
\newtheorem{remark}{Remark}

\journal{Applied Mathematics Letters}

\begin{document}

\begin{frontmatter}



\title{Mixed-Dimensional Geometric Coupling\\ of Port-Hamiltonian Systems}


\author[AMNA]{Jens Jäschke\corref{cor1}}
\ead{jaeschke@uni-wuppertal.de}
\author[FA]{Nathanael Skrepek}
\author[AMNA]{Matthias Ehrhardt}

\cortext[cor1]{Corresponding author}

\address[AMNA]{Bergische Universität Wuppertal, Applied Mathematics and Numerical Analysis, Gaußstraße 20, 42119 Wuppertal, Germany}

\address[FA]{Bergische Universität Wuppertal, Functional Analysis, Gaußstraße 20, 42119 Wuppertal, Germany}



\begin{abstract}
We propose a new interconnection relation for infinite-dimensional port-\linebreak{}Hamiltonian systems that enables the coupling of ports with different spatial dimensions by integrating over the the surplus dimensions. 
To show the practical relevance, we apply this interconnection to a model system of an actively cooled gas turbine blade.
We also show that this interconnection relation behaves well with respect to a discretization in finite element space, ensuring its usability for practical applications.
\end{abstract}




\end{frontmatter}



\section{Motivation}
Scientific models are inherently approximations of reality, and removing unnecessary details can greatly simplify the resulting model.
These simplifications often involve reducing the spatial dimensions of the model: 
A fluid flowing through a pipe is often modelled in 1D rather than using the full 3D Navier-Stokes equations. 
Electronic components such as capacitors and resistors are commonly modelled as 0D elements.
When the interfaces of the subsystems have the same dimension, there are formalisms such as \textit{Port-Hamiltonian Systems} (PHS) that treat the interconnection of these systems in a fairly general way.

However, it becomes difficult when the subsystems have different spatial dimensions. 
For example, modeling a one-dimensional pipe flow that interacts with its environment via the pipe walls requires coupling a 1D interface (the fluid flow) with a 2D interface (the pipe walls). 
Coupling the pipe walls to a lumped-parameter model for the temperature of the room in which they are located requires coupling the 2D pipe surface to a zero-dimensional system.

In the following sections, we will attempt to formulate an energy-conserving connection of two port-Hamiltonian systems where the connected ports do not have the same spatial dimension.

\section{Motivating Example: Cooled Gas Turbine Blade}
\label{sec:motivating_example}
\begin{figure}
    \centering
    \begin{subfigure}[b]{0.45\textwidth}
    \centering
    \includegraphics[width=\textwidth]{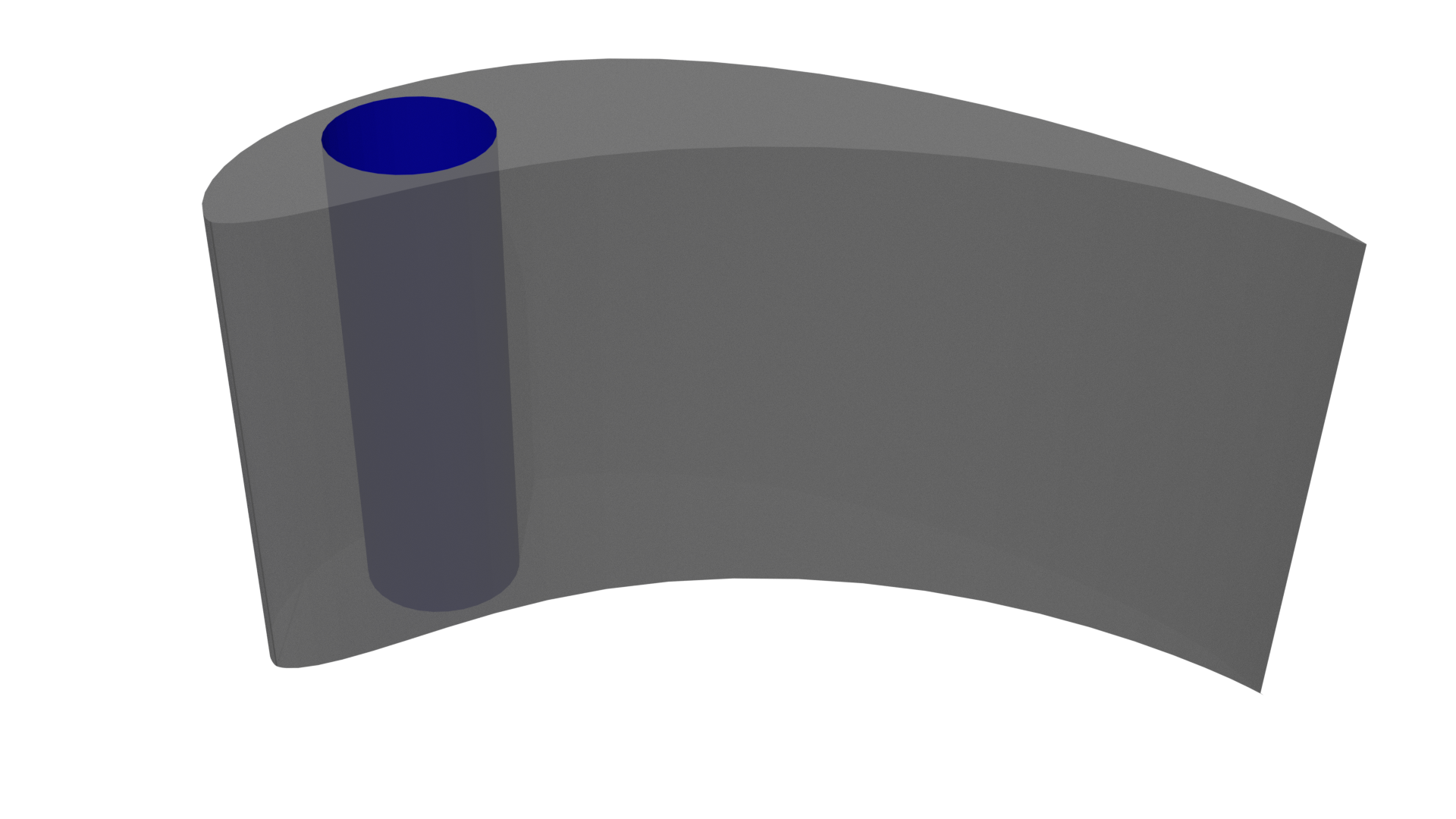}
    \caption{3D view}
    \label{fig:turbine_blade_model_3d}
    \end{subfigure}
    \begin{subfigure}[b]{0.45\textwidth}
    \centering
    \includegraphics[width=\textwidth]{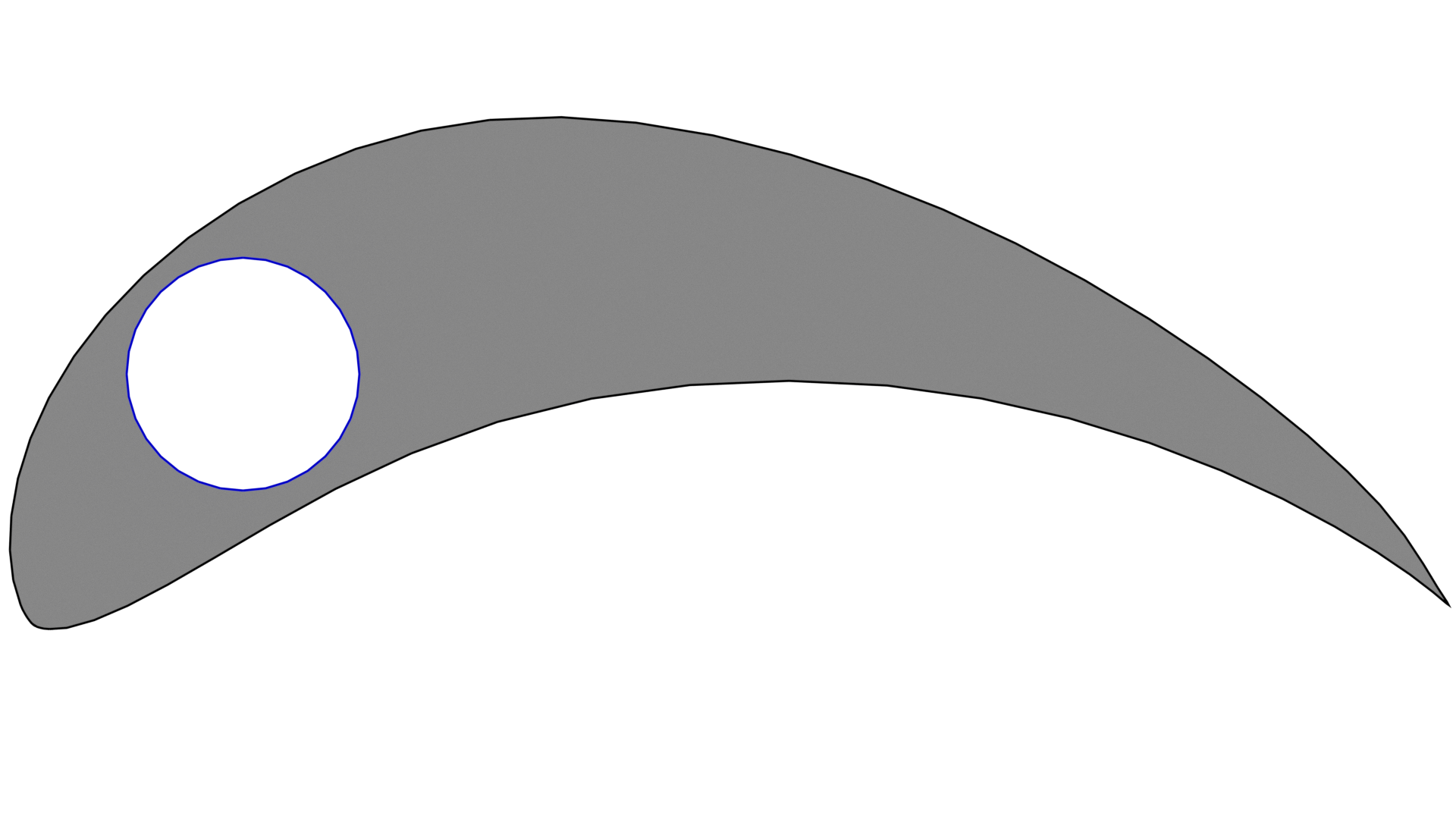}
    \caption{Top view}
    \label{fig:turbine_blade_model_top}
    \end{subfigure}
    \caption{Simple model of a cooled turbine blade, with the cooling channel in blue.}
    \label{fig:turbine_blade_model}
\end{figure}

Consider the heat flow in a gas turbine blade cooled by an internal cooling channel, as shown in \Cref{fig:turbine_blade_model}.
We can model this as two interconnected subsystems: the heat conduction within the metal of the turbine blade and the coolant flow within the cooling channel.
For more information and a discussion of a greatly simplified version of this system, see \cite{jaschke2021preprint}.

Heat conduction in the metal is, of course, modelled by a heat equation:
\begin{equation*}
   \rho c \pp{T}{t}(x,t) = \divergence \bigl( \lambda \grad T(x,t)\bigr)
\end{equation*}
The formulation as a port-Hamiltonian system closely follows \cite{Serhani2019}, choosing the thermal energy $Q$ as Hamiltonian
\begin{equation*}
    Q(t) = \int_\Omega q\bigl(s(x,t)\bigr)\dif{x},
\end{equation*}
and considering the thermal energy density $q$ as a function of the entropy density $s$ such that the thermodynamic relation $\delta_s Q = \dd{q}{s} = T$ is satisfied.
Taking $s$ as a state variable, we obtain the usual flow $f_s = \pp{s}{t}$ and the corresponding effort $e_s = T$. 
As additional flows and efforts we choose the entropy flux $e_\Phi = \Phi_S$, as well as $f_\Phi = -\grad(T)$, $f_\sigma = T$ and $e_\sigma= -\grad(\frac{1}{T})\Phi_Q$ with the heat flux $\Phi_Q$ to obtain the port Hamiltonian system
\begin{align*}
    \begin{pmatrix}
        f_s \\ f_\Phi \\ f_\sigma
    \end{pmatrix}
    =
    \begin{pmatrix}
        0 & -\divergence & -1 \\
        -\grad & 0 & 0 \\
        1 & 0 & 0
    \end{pmatrix}
    \begin{pmatrix}
        e_s \\ e_\Phi \\ e_\sigma
    \end{pmatrix}.
\end{align*}
Since this has two algebraic equations, we add the two \textit{closure relations}
\begin{align*}
    e_s e_\Phi = \lambda f_\Phi
    \quad\text{and}\quad
    f_\Phi e_\Phi = - f_\sigma e_\sigma,
\end{align*}
the former being Fourier's law and the latter expressing the relation between heat flux $\Phi_Q$ and entropy flux $\Phi_S$. 
Finally, we choose the input $u = T|_{\partial \Omega}$ and output $v = -\left(\Phi_S \vec{n}\right)|_{\partial \Omega}$ with $\vec{n}$ being the surface normal vector.

The coolant flow in the cooling channel is modelled as a 1D compressible fluid.
This is consistent with common practice in engineering, since cooling channels in practice are small, irregularly shaped, and exhibit highly turbulent flow, making full 3D flow models infeasible for practical applications and requiring the use of 1D parameter models, such as those presented in \cite{meitner1990computer}.
A 1D model also allows us to use the formulation of irreversible PHS with boundary control presented in \cite{ramirez2022}.
We choose the specific volume $\varphi = \frac{1}{\rho}$, the velocity $v$ and the entropy density $s$ as state variables, and the Hamiltonian 
\begin{equation*}
    H(v,\varphi,s) = \int_a^b \Bigl(\frac{1}{2}v^2 + u(\varphi,s) \Bigr)\dif{z},
\end{equation*}
where the internal energy density $u$ fulfils the Gibbs relation $\dif{u} = -p\dif{\varphi} + T\dif{s}$.
We can then formulate the quasi-Hamiltonian system
\begin{align*}
    \begin{pmatrix}
        \pp{\varphi}{t} \\ \pp{v}{t} \\ \pp{s}{t}
    \end{pmatrix}
    &=
    \begin{pmatrix}
    0 & \pp{\cdot}{z} & 0 \\
    \pp{\cdot}{z} & 0 & -\frac{fv}{T} \\
    0 & \frac{fv}{T} & 0
    \end{pmatrix}
    \begin{pmatrix}
    -p \\ v \\ T
    \end{pmatrix}
    +
    \begin{pmatrix}
    0 \\ 0 \\ 1
    \end{pmatrix}
    w(z,t),\\
    y &= 
    \begin{pmatrix} 
    0 & 0 & 1  
    \end{pmatrix}
    \begin{pmatrix}
    -p \\ v \\ T
    \end{pmatrix} = T,
\end{align*}
with the appropriate boundary conditions.
This system is an infinite-dimensional irreversible port Hamiltonian system as defined in \cite[definition~1]{ramirez2022}.

Coupling the two systems using the usual interconnections for PHS does not work because the spatial dimensions do not match:
The boundary port of the heat equation is 2D, while the distributed port of the cooling channel is only 1D.
We need a new interconnection to compensate for this dimensional mismatch.

\section{Proposition: Mixed-Dimensional Geometric Coupling}\label{sec:coupling}

\begin{definition}[Dirac structure \cite{leGorrec2005}]
Let $\mathcal{F}$ be a linear space, $\mathcal{E}$ its dual and
    $\scp{ \cdot }{ \cdot }\colon \mathcal{E}\times\mathcal{F} \to \Reals$ their dual product. 
    Further let
    \begin{equation*}
    \sdp*{%
    \begin{pmatrix}
    e_1\\f_1
    \end{pmatrix}}
    {\begin{pmatrix}
    e_2\\f_2
    \end{pmatrix}}
    = \scp{e_1}{f_2} + \scp{e_2}{f_1}
    \quad
    \begin{pmatrix}
    e_1\\f_1
    \end{pmatrix},
    \begin{pmatrix}
    e_2\\f_2
    \end{pmatrix}
    \in \mathcal{E}\times\mathcal{F}.
    \end{equation*}
    Then $\mathcal{D} \subset (\mathcal{E} \times \mathcal{F})$ is a \emph{Dirac structure} if $\mathcal{D} = \mathcal{D}^\bot$ with
    \begin{equation*}
       \mathcal{D}^\bot = \{a \in \mathcal{E}\times \mathcal{F} \;|\; \sdp{a}{b} = 0 \quad \forall\ b \in \mathcal{D}\}.
    \end{equation*}
\end{definition}

\begin{theorem}\label{thm:coupling}
Let $\Gamma_1 \subseteq \Reals^n$ compact, $\Gamma_2 \subset \Reals^m$ compact and $\Gamma \coloneqq \Gamma_1 \times \Gamma_2 \subseteq \Reals^{n+m}$.
Further let $\mathcal{F}=\mathrm{L}^2(\Gamma_1)\times\mathrm{L}^2(\Gamma)$ and $\mathcal{E}=\mathcal{F}^*$ its dual. 
Note that we have for $x\in\Gamma$ the decomposition $x = (x_1,x_2)$ with $x_1\in\Gamma_1$ and $x_2\in\Gamma_2$.
Finally, let 
\begin{equation*} \label{eq:integral_operator}
    A \colon\left\{
    \begin{array}{rcl}
        \mathrm{L}^2(\Gamma) & \to & \mathrm{L}^2(\Gamma_1), \\
        u & \mapsto & \int_{\Gamma_2} u(\cdot,x_2) \,\dif{x_2},
    \end{array}
    \right.
\end{equation*} 
and the embedding
\begin{equation*}
B \colon\left\{
\begin{array}{rcl}
    \mathrm{L}^2(\Gamma_1) &\to & \mathrm{L}^2(\Gamma),\\
    v &\mapsto & v.
\end{array}
\right.
\end{equation*}
The previous operator has to be understood as $(Bv)(x_1,x_2) = v(x_1)$.
Then 
\begin{equation*}
J \colon\left\{
\begin{array}{rcl}
    \mathcal{E} &\to & \mathcal{F},\\
    e &\mapsto & 
    \begin{pmatrix}
    0 & -A \\ B & 0
    \end{pmatrix}
    \begin{pmatrix}e_1\\e_2\end{pmatrix},
\end{array}
\right.
\end{equation*}
induces a Dirac structure
\begin{equation*}
    \mathcal{D} = \bigl\{ (e,f) \in \mathcal{E} \times \mathcal{F} \;|\; f=Je \bigr\}.
\end{equation*}
\end{theorem}

Note that $u \in \mathrm{L}^{2}(\Gamma)$ implies $u(\cdot,x_2) \in \mathrm{L}^{2}(\Gamma_{1})$ for almost every $x_{2} \in \Gamma_{2}$. 
Moreover, by the triangle inequality and Cauchy-Schwarz inequality
\begin{align*}
    \norm{Au}^2_{\mathrm{L}^2(\Gamma_1)}
    &=\int_{\Gamma_1} \abs*{\int_{\Gamma_2} u(x_1,x_2) \dif{x_2}}^{2} \dif{x_{1}}
    \leq \int_{\Gamma_1} \Big(\int_{\Gamma_2} 1\cdot \abs{u(x_1,x_2)}  \dif{x_2}\Big)^{2} \dif{x_{1}} \\
    &\stackrel{C.S.}{\leq} \abs{\Gamma_2}\int_{\Gamma_1} \int_{\Gamma_2} \abs{u(x_1,x_2)}^{2} \dif{x_2} \dif{x_{1}} = \abs{\Gamma_{2}} \norm{u}_{\mathrm{L}^{2}(\Gamma)}^2,
\end{align*}
where $\abs{\Gamma_{2}}$ is the measure of $\Gamma_{2}$. 
Hence, the operator $A$ is well-defined. Note that this holds true for any finite measure on $\Gamma_2$. In particular we will later use surface measures.

\begin{proof}
    Determine the adjoint operator of $B$: For $f \in \mathrm{L}^2(\Gamma)$, $v \in \mathrm{L}^{2}(\Gamma_{1})$ we have
    \begin{equation*}
        \begin{split}
            \scp{f}{Bv}_{L^2(\Gamma)} 
            &= \int_{\Gamma_1} \int_{\Gamma_2} f(x_1,x_2) v(x_1) \dif{x_2} \dif{x_1} 
            = \int_{\Gamma_1}\Bigl( \int_{\Gamma_2} f(x_1,x_2) \dif{x_2} \Bigr) v(x_1)\dif{x_1} \\
            &= \scp*{\int_{\Gamma_2} f(\cdot,x_2) \dif{x_2}}{v}_{L^2(\Gamma_1)} 
            = \scp*{B^* f}{v}_{L^2(\Gamma_1)} 
            = \scp*{A f}{v}_{L^2(\Gamma_1)}.
        \end{split}
    \end{equation*}
    Since $A=B^*$ holds, $J$ is skew-symmetric and $\mathcal{D}$ is a Dirac structure \cite{vanderSchaft2002}. 
\end{proof}

\section{Coupled Example System}\label{sec:coupling_example}


To apply the coupling described in \Cref{sec:coupling} to the system of \Cref{sec:motivating_example},
we first split the boundary of the heat equation domain $\partial \Omega = \Gamma$ into 
an external part $\Gamma_{\rm ext}$, which connects to the outside of the blade and is disregarded here, and an internal part $\Gamma_{\rm int}$ which denotes the wall of the cooling channel and will be coupled to the coolant flow.

As the cooling channel is modelled as a tube, it can obviously be decomposed into $\Gamma_{\rm int} = \Gamma_1 \times \Gamma_2$ as in
\Cref{thm:coupling},
with $\Gamma_1$ containing the axial coordinate (along the flow direction) and $\Gamma_2$ the azimuthal coordinate, i.e.\ describing the circumference. 
Now we can choose the following interconnection
\begin{align}
    u = T_h(x) = T_c(z) = y \quad \text{and} \quad
    w 
    = \int_{\Gamma_2} \Phi_S(x)\, \vec{n}\, \dif{x_2} 
    = - \int_{\Gamma_2} v \,\dif{x_2},
\label{eq:interconnection_relations}
\end{align}
with $x \in \Gamma_{\rm int}$ and $\scp{x}{\hat{z}} = z$, where $\hat{z}$ is the unit vector in $z$-direction.
This interconnection has exactly the form given in
\Cref{thm:coupling}.
Since it is an energy preserving interconnection, the coupled system is a (quasi-)Hamiltonian system and would be a port-Hamiltonian system if both sub-systems were PHS.

We note that this interconnection is also physically meaningful. 
The temperature $T$, an intensive quantity, of the points that are in contact with each other is the same, while the entropy flux $\Phi_S$, an extensive quantity, is integrated and has the expected sign change.

\section{Finite Element Discretization}

The interconnection proposed in \Cref{sec:coupling}
can be easily discretized with a finite element discretization. The result will then be a finite-dimensional Dirac structure, as we will see in this section. 

Let us assume that we have finite element discretizations for both sub-systems, with ${\psi_i}$ the basis functions on the boundary of the higher-dimensional system (the heat equation in our example), and ${\chi_i}$ the basis functions of the lower-dimensional system (the compressible cooling fluid in our example).
We can then approximate the input $u$ and output $v$ of the first system, and the input $w$ and output $y$ of the second system as
\begin{gather*}
\label{eq:fe_approximations}
\begin{aligned}
    u &\approx \sum_i \psi_i(x) u_i(t)  = \Psi^\top(x) \underline{u}(t), & v &\approx \sum_i \psi_i(x) v_i(t) = \Psi^\top(x) \underline{v}(t), \\
    w &\approx  \sum_i \chi_i(x_1) w_i(t) = X^\top(x_1) \underline{w}(t), & y &\approx \sum_i \chi_i(x_1) y_i(t) = X^\top(x_1) \underline{y}(t).
\end{aligned}
\end{gather*}
Remembering that $x = (x_1, x_2)^\top$ and applying these approximations to the continuous interconnection relations of \Cref{eq:interconnection_relations}
results in
\begin{align}\label{eq:approx_interconnection}
    X^\top(x_1) \underline{w}(t) = -\int_{\Omega_2} \Psi^\top(x) \underline{v}(t) \dif{x_2} , \quad \text{and} \quad
    \Psi^\top(x) \underline{u}(t) = X^\top(x_1) \underline{y}(t).
\end{align}
We now take the weak form of \Cref{eq:approx_interconnection} to obtain the discretized forms of the interconnection relations
\begin{gather*}
\begin{aligned}
    M_\chi \underline{w}(t) &= \int_{\Gamma_1} X(x_1) X^\top(x_1) \underline{w}(t) \dif{x_1} = - \int_{\Gamma_1} X(x_1) \int_{\Gamma_2} \Psi^\top(x) \underline{v}(t) \dif{x_2} \dif{x_1} \\
    &= -\int_{\Gamma_1} X(x_1) \widehat{\Psi}^\top(x_1) \underline{v}(t) \dif{x_1} = - D_\chi \underline{v}(t)
\end{aligned}
\end{gather*}
and
\begin{gather*}
\begin{aligned}
       M_\psi \underline{u}(t) &= \int_\Gamma \Psi(x) \Psi^\top(x) \underline{u}(t) \dif{x} = \int_\Gamma \Psi(x) X^\top(x_1) \underline{y}(t) \dif{x}\\
        &= \int_{\Gamma_1} \left( \int_{\Gamma_2} \Psi(x) \dif{x_2} \right) X^\top(x_1) \underline{y}(t) \dif{x_1} \\
        &= \int_{\Gamma_1} \widehat{\Psi}(x_1) X^\top(x_1) \underline{y}(t) \dif{x_1} = D_\psi \underline{y}(t).
\end{aligned}
\end{gather*}
Since $D_\psi = D_\chi^\top$, the discretized interconnection relation
\begin{align*}
    \begin{pmatrix}
    M_\chi & 0 \\
    0 & M_\psi
    \end{pmatrix} 
    \begin{pmatrix}
    \underline{u}(t) \\ \underline{w}(t)
    \end{pmatrix}
    &=
    \begin{pmatrix}
    0 & -D_\chi \\
    D_\psi & 0
    \end{pmatrix}
    \begin{pmatrix}
    \underline{v}(t) \\ \underline{y}(t)
    \end{pmatrix}
\end{align*}
is a Dirac structure.

\begin{remark}
The integration over $\Gamma_2$ will not expand the support of the basis functions $\widehat{\psi}_i$ in $x_1$-direction. Therefore, the matrix $D_\chi$ will still be sparse, although less sparse than the matrix $M_\chi$.
\end{remark}

\section{Conclusion}
It is possible to couple port-Hamiltonian systems of different spatial dimensions if the interconnecting ports do not have the same spatial dimension. 
The proposed interconnection structure forms a Dirac structure and thus ensures that the resulting overall system again forms a port-Hamiltonian system. 

Application to an example system has shown that the interconnection has practical use and a physically meaningful interpretation when the ports consist of both extensive and intensive variables.
This is usually the case for physically motivated 
port-Hamiltonian systems, but cannot be guaranteed in general. 

Finally, we showed that the interconnection behaves well with respect to the discretization in finite element space, leading to a finite-dimensional Dirac structure.



\bibliographystyle{unsrtnat}
\bibliography{bibliography}
 
\end{document}